 \documentclass[12pt]{amsart}

 \usepackage{mathrsfs,amsfonts}
 \usepackage{color}

 \newtheorem{theorem}{Theorem}[section]
 \newtheorem{lemma}[theorem]{Lemma}
 \newtheorem{corol}[theorem]{Corollary}
 \newtheorem{prop}[theorem]{Proposition}
 \newtheorem{example}[theorem]{Example}
 
 \newtheorem{con1}[theorem]{Condition}
 \newtheorem{remark}[theorem]{Remark}

 \def\blemma{\begin{lemma}\sl{}\def\elemma{\end{lemma}}}
 \def\bproposition{\begin{prop}\sl{}\def\eproposition{\end{prop}}}
 \def\btheorem{\begin{theorem}\sl{}\def\etheorem{\end{theorem}}}
 \def\bcorollary{\begin{corol}\sl{}\def\ecorollary{\end{corol}}}

 \def\beqlb{\begin{eqnarray}}\def\eeqlb{\end{eqnarray}}
 \def\beqnn{\begin{eqnarray*}}\def\eeqnn{\end{eqnarray*}}

 \def\proof{\noindent{\it Proof.~~}}\def\qed{\hfill$\Box$\medskip}

 \def\<{\langle}\def\>{\rangle}

 \def\mcr{\mathscr}\def\mbb{\mathbb}

 \def\ar{\!\!&}

\begin{document}

\

\bigskip\bigskip

\centerline{\Large\bf Integral functionals for spectrally positive }
\

\centerline{\Large\bf L\'{e}vy processes }

\bigskip\bigskip

\centerline{Pei-Sen Li $^{a)}$ and Xiaowen Zhou $^{b)}$}

\smallskip
\centerline{Institute for Mathematical Sciences, Renmin University of China,}

\centerline{59 Zhongguancun Street, Beijing, P. R. China}

\centerline{and}

\centerline{Department of Mathematics and Statistics, Concordia University,}

\centerline{1455 De Maisonneuve Blvd. W., Montreal, Canada}

\centerline{\footnotesize{  a) peisenli@mail.bnu.edu.cn, b) xiaowen.zhou@concordia.ca}}
\bigskip


\bigskip

\centerline{Abstract}
\bigskip

 We find necessary and sufficient conditions for almost sure  finiteness of  integral functionals  of spectrally positive L\'{e}vy processes. Via Lamperti type transforms, these results can be applied to obtain new integral tests on extinction and explosion behaviors for a class of  continuous-state nonlinear branching processes.

\smallskip



\bigskip

\section{Introduction}

\setcounter{equation}{0}
For a diffusion process or a L\'evy process $Z=(Z_t)_{t\ge0}$ taking values in an open interval $I = (l, r)$,
let $f$ be a non-negative, measurable and locally bounded function on $I$, and introduce
\beqlb\label{func}
A_t(f) := \int^t_0 f(Z_s ) ds.
\eeqlb
Denote by  $\zeta:=\inf\{t>0: Z_t> r~\mbox{or}~Z_t< l\}$ the exiting time of $Z.$
The ultimate value of this additive functional, i.e., $A_\zeta(f)$, is often called a \emph{perpetual integral functional}. In this short note we are interested in finding necessary and sufficient conditions for the a.s. finiteness of $A_\zeta (f)$ for certain choices of $\zeta$.

The perpetual integral functionals appear naturally in various applications. One of the common issues in studying Markov processes is the question that whether a random time-change of the process is again conservative. The time-changed process is usually realized as the generalized inverse of a positive additive functional. As a consequence, the time-changed process is conservative if and only if $A_\zeta (f)=\infty$ a.s.  In insurance
mathematics,  the functional $A_\zeta(f)$ may be interpreted, in a suitable model, as the present value
of a continuous stream of perpetuities, see Dufresne \cite{Du90}.
In the stochastic models of ecology, the functional $A_\zeta(f)$ can  be regarded as the total population of a
species.

The finiteness of (\ref{func}) has  been studied by various authors. The well known  \emph{Engelbert-Schmidt zero-one law} (see \cite[Proposition 3.6.27]{KaS91}) states that for a Brownian motion $(B_t)_{t\ge 0}$ and any nonnegative Borel function $f$ the following three statements are equivalent:
\begin{itemize}
\item[(i)]
 $$\mbb{P}\left(\int^t_0 f(B_s)ds<\infty ~\mbox{for all}~t\in[0, \infty)\right)> 0,$$
\item[(ii)]
 $$\mbb{P}\left(\int^t_0 f(B_s)ds<\infty ~\mbox{for all}~t\in[0, \infty)\right)=1,$$
\item[(iii)]
the function $f$ is locally integrable on $\mbb{R}$.
\end{itemize}
This important property has a plenty of applications. For example, it constitutes an important step in the Engelbert-Schmidt construction of weak solutions of one-dimensional SDEs.
When $Y$ is a Brownian motion with positive drift,  it is known that $A_\infty(f)$ is finite a.s. if and only
if $f$ is integrable at $\infty$ (see Engelbert and Senf [9] and Salminen and Yor [18]).
When $Z$ is a diffusion process on an interval $I=(r, l)$. The finiteness of the perpetual integral functional of $Z$, i.e. $A_\zeta(f)$, is considered in  Khoshnevisan et al. \cite{Kh06} and where the proofs are based on  Khasminskiis  theorem.

When $Z$ is a L\'evy process, the fact that $Z$ may have jumps makes things more involved. Suppose that $f$ is a positive locally integrable function and $Z$ is a L\'{e}vy process such that $\mbb{E}[Z_1]\in(0, \infty)$ and its local time exists. D\"{o}ring and Kyprianou \cite{Dor16} find the following condition
\beqnn
\int^\infty_0 f(Z_s)ds<\infty \,\, a.s. \Leftrightarrow \int^{\infty} f(x)dx<\infty.
\eeqnn
Their proof is based on Jeulin's lemma.
Recently, Kolb and Savov \cite{KoS19} improve the above result by removing the assumption of the existence of local time. They establish that
\beqnn
\int^\infty_0 f(Z_s)ds<\infty \,\, a.s. \Leftrightarrow \int^{\infty} f(x)U(dx)<\infty,
\eeqnn
where $U$ is the potential measure of $Z$.

When the function $f$ is an exponential function, i.e. $f(x) = e^{-\theta x}$ for some $ \theta> 0$, the exponential integral functional has drawn the attention of many authors in recent years; see, e.g.\cite{EricMaller, KunPar, Kon, Pardo09, Pardo13}. The exponential integral functional plays an important role in many domains, e.g., mathematical finance, branching process in random environments, tail probability, especially in the case of L\'evy process with one-sided jumps (see for example \cite{Patie12, Patie15}). We refer to Bertoin and Yor \cite{BerYor05}  for a survey of this area.

However, in the existing literature, when $Z$ is a L\'evy process, only functionals with deterministic upper limit have been considered so far. To the best of our knowledge, there has been no systematic discussion on the a.s. finiteness of functional with a stopping time as the upper limit for integral.
In this work, we are going to fill this gap. For a spectrally positive L\'evy process $Z$, i.e.  a one-dimensional L\'evy process $Z$ with no negative jumps, we obtain an explicit integral test on the finiteness of $A_\zeta(f)$, where $\zeta$ denotes the exiting time of $(0, \infty)$, i.e. $\zeta=\inf\{t>0: Z_t=0\}$ with the convention $\inf\emptyset=\infty$.

The assumption of non-negative jump allows us to use the tools and technics available for one-sided L\'evy processes,  such as scale functions. On the other hand, the spectrally positive L\'evy process is of particular interest since in this case, the  processes obtained by time-changing process $Z$ forms a class of important processes with rich mathematical structures; see for example the continuous-state nonlinear branching processes. We will discuss the time-changed processes in details in the next section.

The remainder of this paper is organized as follows. In Section 2 we introduce the necessary notations
 and present the main results on the Engelbert-Schmidt type zero-one law for spectrally positive L\'evy process.   We also apply the main results to find  integral tests on the boundary behaviors of the time-changed processes. The last section is devoted to the proofs of the main results.

\section{Notations and main results}

In this section we first recall the L\'evy-It\^{o} decomposition of spectrally positive L\'evy process and some necessary notations.
Let $b$ and $c\geq 0$ be constants and $\pi(dz)$ be a $\sigma$-finite measure on $(0,\infty)$ satisfying
 \beqnn
\int_{(0,\infty)} (1\wedge z^2) \pi(dz)< \infty.
 \eeqnn
  Given a Brownian motion $(B_t)_{ t\ge 0}$   and an independent Poisson random measure $N(ds,dz)$ on $(0,\infty)\times (0,\infty)$ with intensity $ds\pi(dz)$, both defined on a filtered probability space $(\mcr{F}, (\mcr{F}_t)_{t\geq 0},  \mbb{P})$,  a spectrally positive L\'{e}vy process $(Z_t)_{ t\ge 0}$ started at $x$ can be represented as
 \beqlb\label{0.3}
Z_t = x - bt + \sqrt{2c} B_t + \int^t_0 \int_{(0,1]} z \tilde{N}(ds,dz)
+\int^t_0 \int_{(1,\infty)} z N(ds,dz),
 \eeqlb
where $\tilde{N}(ds,dz):= N(ds,dz)-ds\pi(dz)$ denotes the compensated Poisson random measure on $(0,\infty)\times (0,1] $. {Throughout the paper we assume that process $Z$ is not a subordinator.}

 Write $\mbb{P}_x(\cdot):=\mbb{P}(\cdot|Z_0=x)$ and $\mbb{E}_x$ for the corresponding expectation.
The Laplace exponent $\psi : [0, \infty) \rightarrow (-\infty, \infty)$ of the spectrally negative L\'evy process $-Z$ is specified by
\[\mbb{E}_x[e^{-\lambda Z_t}] = e^{-\lambda x+\psi(\lambda)t},\qquad \lambda,t\geq 0,\]
where $\psi$ has the expression of
\beqlb\label{0.01}
\psi(\lambda) = b\lambda+c\lambda^2+\int_{(0,\infty)}(e^{-\lambda u} - 1 +
\lambda u 1_{(0,1]}(u)) \pi(du),\quad \lambda\ge 0.
 \eeqlb
{By the assumption that $Z$ is not a subordinator,} we see that there exists a $\lambda\in (0,\infty)$ such that $\psi(\lambda)> 0$ 
In that case, it is known that
the Laplace exponent $\psi$ is strictly convex and tends to $\infty$ as $\lambda\rightarrow\infty$.
Let
\beqnn
\Phi(0):=\inf\{\lambda > 0: \psi(\lambda) > 0\}<\infty.
\eeqnn

Recall that we restrict $Z$ on the positive half line. Then the life time of $Z$ is just $\zeta$, the first hitting time of $0$.
By \cite[Theorem 3.12]{Kyp14} we have
\beqlb\label{a00}
\mbb{P}_x(\zeta<\infty)=e^{-\Phi(0)x}\qquad\mbox{and}\qquad \mbb{P}_x(\zeta=\infty)=1-e^{-\Phi(0)x}.
\eeqlb
Here, we stress that $\mbb{P}_x(\zeta<\infty)$ is always strictly positive in our setting. On the other hand, we have $\mbb{P}_x(\zeta=\infty)>0$ if and only if $\Phi(0)>0.$

We are
concerned with the almost sure finiteness of $A_\zeta(f)$ under conditional probabilities $\mbb{P}_x(
\cdot|\zeta<\infty)$ and $\mbb{P}_x(\cdot|\zeta=\infty)$, respectively.
To illustrate our main contributions, we present the following statements for the
 Engelbert-Schmidt type zero-one law for spectrally positive L\'evy processes. The proofs are  deferred to Section 3.

\btheorem\label{t1.1}
 Let $f$ be a strictly positive function on $(0,\infty)$ satisfying $\sup_{x\geq \varepsilon}f(x)<\infty$ for any $\varepsilon>0$. Then for any $x>0$,
the following four statements are equivalent:
\begin{itemize}
\item[(i)]
$$
\mbb{P}_x(A_{\zeta}(f)<\infty|\zeta<\infty)> 0;
$$
\item[(ii)]
$$
\mbb{P}_x(A_{\zeta}(f)<\infty|\zeta<\infty)=1;
$$
\item[(iii)]
$$
\mbb{E}_x\Big[\int^{\zeta}_0 f(Z_{t})e^{-\lambda Z_{t}}dt|\zeta<\infty\Big]<\infty \qquad\mbox{for some/all}~\lambda>0;
$$
\item[(iv)]
\beqnn
\int^{\infty}\frac{f(1/\lambda)}{\lambda\psi(\lambda)}d\lambda<\infty.
\eeqnn
\end{itemize}
\etheorem

 For any $y\ge0$, set $\tau^-_y:=\inf\{t\ge0: Z_t\le y\}$ with the convention $\inf\emptyset=\infty$. Note that $\zeta=\tau^-_0$ by definition.

\btheorem\label{t1.2}
Suppose that $\Phi(0)>0$ and $f$ is a strictly positive and decreasing function on $(0,\infty)$. Then the following three statements are equivalent:

(i) $$
\mbb{P}_x(A_{\infty}(f)<\infty|\zeta= \infty)> 0\qquad\mbox{for all} ~x>0;
$$

(ii) $$
\mbb{P}_x(A_{\infty}(f)<\infty|\zeta= \infty)=1\qquad\mbox{for all} ~x>0;
$$

(iii)
 $$
 \mbb{E}_x[A_{\tau^-_y}(f)]<\infty \qquad\mbox{for all}~x>0~\mbox{and all}   ~ y\in (0,x).
$$
Moreover,  if we further assume that $\psi'(0)<\infty$, then the above three statements are all equivalent to

(iv)
 $$
\int^\infty f(y)dy<\infty.
$$
\etheorem
In the above theorem, the assumption $\psi'(0)<\infty$ could be removed at the price of
 an additional requirement on $f$, and the integral test becomes a little more complicated.
\bcorollary\label{cor2xx}
Suppose that $\Phi(0)>0$ and function $f$  can be represented as the Laplace transform of a nonnegative function $g$, i.e. $f(x)=\int^\infty_0 e^{-xz}g(z)dz$ for $x>0$. Then the  statements (i), (ii) and (iii) in Theorem \ref{t1.2} are all equivalent to
\[
\int_{0+}\frac{g(\lambda)}{\psi(\lambda)}d\lambda>-\infty.
\]
\ecorollary

The integral functionals are very useful in the study of explosion and extinction behaviors of the following SDE, whose solution can be represented as a  time changed spectrally positive L\'{e}vy process. Consider a filtered probability space $(\Omega, \mathscr{G}, \mathscr{G}_t, \mbb{P})$ satisfying the usual hypotheses. Let $\{W_t\}_{t\ge0}$  be an  $(\mathscr{G}_t)$-Brownian motion. Let $M(ds,dz,du)$ be an independent $(\mcr{G}_t)$-Poisson random measure on $(0,\infty)\times (0,\infty)\times (0,\infty)$ with intensity $ds\pi(dz)du$ and let $\tilde{M}(ds,dz,du)$ be the corresponding compensated Poisson random measure. We consider the nonnegative solution of the stochastic integral equation
\beqlb\label{e0.1}
X_t\ar=\ar x+\sqrt{2c}\int^{t}_0 \sqrt{1/f(X_s)} dW_s + \int_0^t\int_{(0,1]} \int^{1/f(X_{s-})}_0 z \tilde{M}(ds,dz,du)\cr
 \ar\ar\quad
-\, b \int^{t}_0 1/f(X_s) ds + \int_0^{t}\int_{(1,\infty]}\int^{1/f(X_{s-})}_0 z M(ds,dz,du).
 \eeqlb
 By  a solution $X= (X_t)_{t\geq 0}$ to equation (\ref{e0.1}) we mean a c\`adl\`ag $(0,\infty)$-valued $(\mcr{F}_t)$-adapted process satisfying~(\ref{e0.1}) up to times
$$\zeta_n := \inf\{t\geq0: X_t\ge n \,\,\text{ or} \,\, X_t\le 1/n\}$$
 for all $n>0$ and  both   $0$ and $\infty$ are absorbing boundaries for $X$. Applying Ikeda and Watanabe \cite[Theorem 9.1]{IkW89} we see that, if $1/f$ is locally Lipschitz, then  SDE {\rm(\ref{e0.1})} has a pathwise unique solution.  Note that if $(Z_t)_{t\ge0}$ is a spectrally one-sided $\alpha$-stable L\'{e}vy process with $\alpha\in(1,2)$, then SDE \eqref{e0.1} can be transformed
into the following form:
 \beqlb\label{e0.2}
X_t\ar=\ar x+\int_0^t \sigma(X_{s-})dZ_t,
 \eeqlb
where $\sigma(X_{s-})=1/f(X_{s-})^{1/\alpha}$. In D\"{o}ring and Kyprianou \cite{Dor18}, they study the boundary behavior of the solution to (\ref{e0.2}), and in their setting, $(Z_t)_{t\ge0}$ is a $\alpha$-stable L\'{e}vy process with $\alpha\in(0,2)$.

The following result on random time change can be proved  using techniques similar to those  in Caballero et al. \cite{CLB09}, where  the function $f(x)=x^{-1}$ is considered.
\bproposition\label{t1.0}
Given any locally bounded and  strictly positive function $f$ on $(0, \infty)$,  for any $t\ge 0$ let
\[{\eta_f(t):= \inf\{s\geq 0: A_{s\wedge\zeta}(f)> t\}    }\]
with the convention $\inf\emptyset=\infty$.
Then $(X_t)_{t\ge0}:= (Z_{\eta_f(t)\wedge\zeta})_{ t\ge 0}$ with convention {$Z_\infty\equiv\infty$ for $\tau^-_0=\infty$} is a weak solution to equation (\ref{e0.1}).
\eproposition

For the case that $f(x)=1/x$, the solution  to (\ref{e0.1}) is called \emph{continuous-state branching process}, which  has attracted the attention of many researchers in the past decades. This process is a model for the evolution of populations. A well known result is that such a process satisfies the \emph{branching property}, i.e.
\[\mbb{E}_{x_1+x_2} [e^{-\lambda X_t}]=\mbb{E}_{x_1} [e^{-\lambda X_t^{(1)}}]\mbb{E}_{x_2} [e^{-\lambda X_t^{(2)}}] \]
 for any $\lambda, t, x_1, x_2>0$, where $X_t^{(1)}$ and $X_t^{(2)}$ denote two independent copies of $X$ with initial values $x_1$ and $x_2$, respectively.
The connection of those processes with spectrally positive L\'evy processes through random time changes is pointed out by Lamperti \cite{Lam67b}.  The continuous-state branching processes also appear in the study of L\'evy trees; see, e.g. Duquesne \cite{Du10}.  A remarkable theory of flows of such processes with applications to flows of Bessel bridges and coalescents with multiple collisions has been developed by Bertoin and Le Gall \cite{BeL00, BeL03, BeL05, BeL06}. The interested readers are referred to Kyprianou \cite{Kyp14}, Li \cite{Li11} and Pardoux \cite{Par16} for reviews of the literature in this subject.

For general locally bounded positive function $f$,  the solution  to (\ref{e0.1})   called \emph{nonlinear branching process} is thus a natural generalization of the continuous-state branching process in which the underlying population evolves non-linearly  governed by the function $1/f$.  Thus, a general nonlinear branching process  no longer satisfies the additive branching property. On the other hand, it could be used to describe the interaction and competition between particles. Such a process allows richer boundary behaviors.  We refer to Li \cite{Li16} and Li et al. \cite{LiYangZhou16} for the criteria of boundary behaviours of this process. The speed of coming down from infinity for such process is studied in Foucart et al. \cite{timechangedsplp} and the speed of explosion is discussed by Li and Zhou \cite{LiB19}.

Thanks to Proposition \ref{t1.0}, we can apply Theorems \ref{t1.1} and \ref{t1.2}  to obtain the following new integral tests on the boundary behaviors for the continuous-state nonlinear branching process $X$ that solves SDE (\ref{e0.1}).
Define the extinction time and the explosion time of $X$ by 
\[T^-_0:=\inf\{t\ge 0: X_t=0\} \quad \mbox{and}\quad  T^+_\infty:=\inf\{t\ge 0: X_t=\infty \},\]
respectively, with the convention $\inf\emptyset=\infty$. We say that process $X$ becomes extinguishing if $T^-_0=\infty$ and $X_t\rightarrow 0$ as $t\rightarrow\infty$. By Proposition \ref{t1.0}, we see that $\eta_f(A_{\zeta}(f))=\zeta$ and
$X_{A_{\zeta}(f)}=Z_{\zeta}$. Then we obtain the following corollary.
\bcorollary\label{Cor}
Given any locally bounded and  strictly positive function $f$ on $(0, \infty)$. We have
$$T^-_0=A_{\zeta}(f)\quad\text{if}\quad  \zeta<\infty\quad\mbox{and}~T^+_\infty=A_{\zeta}(f)\quad\text{ if}\quad  \zeta=\infty.$$
\ecorollary
Thus, the extinction and explosion of the solution $X$  within finite time correspond to the almost sure finiteness of $A_{\zeta}(f)$ conditioned on events $\{\zeta<\infty\}$ and $\{\zeta=\infty\}$, respectively. Similarly, extinguishing occurs for $X$ if $A_{\zeta}(f)=\infty$ given $\{\zeta<\infty\}$. Then, using Theorems \ref{t1.1} and \ref{t1.2}, we can find necessary and sufficient conditions for the process $X$ to die out, to become extinguishing or to explode within finite time.
\bcorollary\label{Cor0}
Given $\theta>0$, let $f(x)=x^{-\theta}$ for $x>0$.

(i) The process $(X_t)_{t\ge 0}$ started at $x>0$ goes extinct in finite time with a positive probability if and only if
\beqnn
\int^{\infty-} \frac{s^{\theta-1}}{\psi(s)} ds <\infty;
\eeqnn
and the process does not go extinct but become extinguishing with a positive probability if  and only if
\beqnn
\int^{\infty-} \frac{s^{\theta-1}}{\psi(s)} ds =\infty.
\eeqnn

(ii) The process $(X_t)_{t\ge 0}$ started at $x>0$ explodes in finite time with a positive probability if and only if $\Phi(0)>0$ and
\beqnn
  \int_{0+}\frac{s^{\theta-1}}{\psi(s)} ds> -\infty.
  \eeqnn
\ecorollary
The above corollary generalizes the results in Grey \cite{Gre74} and Kawazu and in Watanabe \cite{IkW89}, where the  classical continuous-state branching process corresponding to $\theta=1$ is considered.


Proofs of the main results are given in Section 3.

\section{Proofs of the main results}

As the fluctuation theory of the L\'evy process $Z$  plays an important role in our proofs, we recall the definition of the scale function and some of its basic properties. There exists a strictly increasing and positive continuous function $W$ on $[0, \infty)$, called scale function, such that
\begin{equation}\label{defscale}
\int_{0}^{\infty}e^{-\lambda y}W(y)d y=\frac{1}{\psi(\lambda)},\qquad \lambda>\Phi(0).
\end{equation}
Define $W(x)=0$ for $x<0$ and write $W=W^{(0)}$ for short.

Some classical results  for the scale functions for spectrally positive L\'evy processes are summarized in the following lemmas. We refer to Bertoin \cite[Chapter VII]{Ber96}  and Kyprianou \cite[Chapter 8]{Kyp14} for the general theory of scale function.

\begin{lemma}\label{l1}
For $x>0$, the potential measure can be represented as
\beqlb\label{lap1}
\mu_x(dy):=\mbb{E}_x\Big[\int^{\zeta}_0 1_{\{Z_{t}\in dy\}}dt\Big]=[e^{-\Phi(q)x}W(y)-W(y-x)]dy.
\eeqlb
 Moreover,  given  $\varepsilon>0$ there exist $c_1,c_2>0$ such that for any $x\in (0, \varepsilon)$
\begin{equation}\label{boundW}
c_1\frac{1}{x\psi\left(1/x\right)} \leq W(x)\leq c_2\frac{1}{x\psi\left(1/x\right)}.
\end{equation}
\end{lemma}
\proof
We refer to \cite[Theorem 8.1]{Kyp14} for the existence of the scale function and  \cite[Theorem 2.7 (ii)]{Ku11}  for (\ref{lap1}).
In the following, we are going to prove (\ref{boundW}) for completeness.
Let
\[\psi^\natural(\lambda):=\psi(\lambda+\Phi(0)),\,\,\lambda>0\,\,\,\,\text{ and} \,\,\,\, W^\natural(x):=e^{-\Phi(0)x}W(x).\]
Note that, by the definition we have
\beqlb\label{limitW}
\lim_{x\to 0}\frac{W^\natural(x)}{W(x)}=1.
\eeqlb
 From \cite[p.193]{Ber96}  we have that $\psi^\natural$ is again a Laplace exponent and that $W^\natural$ is the corresponding scale function. By the proof of \cite[Propositions VII.10]{Ber96}, we see that there exist $c'_1, c'_2>0$ such that for any $x>0$
  \beqlb\label{boundWn}
c'_1\frac{1}{x\psi^\natural\left(1/x\right)} \leq W^\natural (x)\leq c'_2\frac{1}{x\psi^\natural\left(1/x\right)}.
\eeqlb
Using the following  form of Taylor's formula
 \beqnn
g(y+z)-g(y)-zg'(y)=z^2\int_0^1 g''(y+zv)(1-v)d v\quad z>0,
 \eeqnn
 The Laplace exponent $\psi$ can be expressed as
 $$
 \psi(\lambda)=b\lambda+c\lambda^2+\lambda^2\int_{(0,1]}u^2\pi(du)\int_{0}^{1}e^{-\lambda u}(1-v)dv+\int_{(1, \infty)}(e^{-\lambda u}-1)
 \pi(du),\quad \lambda\ge 0.
 $$
 Then
\beqnn
 \psi^\natural(\lambda)\ar\le\ar b(\lambda+\Phi(0))+c(\lambda+\Phi(0))^2+(\lambda+\Phi(0))^2\int_{(0,1]}u^2\pi(du)\int_{0}^{1}e^{-\lambda u}(1-v)dv\cr
 \ar\ar\qquad\qquad+\int_{(0, \infty)}(e^{-\lambda u}-1)
 \pi(du).
\eeqnn
Applying the monotonicity of $\psi^\natural$, one can check that
\beqlb\label{limit1}
\lim_{\lambda\to\infty}\frac{\psi^\natural(\lambda)}{\psi(\lambda)}= 1.
\eeqlb
Combing (\ref{limitW}), (\ref{boundWn}) and (\ref{limit1})  we prove (\ref{boundW}).
 \qed

Scale functions rarely have explicit expressions.  In the critical stable case, for which $\Psi(\lambda)=\lambda^{\alpha}$ with $\alpha\in (1,2]$, the scale function can be found for instance in \cite[Example 4.17]{Ku11}, and $W(x)=x^{\alpha-1}/\Gamma(\alpha)$, for $x>0$.

\blemma\label{l1'}
For any $x > 0$  and $ \lambda> 0$, we have
\beqlb\label{exp1}
\mbb{E}_x\Big[\int^{\zeta}_0e^{-\lambda Z_t}dt|\zeta<\infty\Big]=\frac{1-e^{-\lambda x}}{\psi(\lambda+\Phi(0))}<\infty,
\eeqlb

%
%
\elemma
\proof
By the Markov property, we have
\beqlb
\mbb{E}_x\Big[\int^{\zeta}_0 e^{-\lambda Z_{t}}dt|\zeta<\infty\Big]
\ar{\color{red}=}\ar
\frac{1}{\mbb{P}_x(\zeta<\infty)}\int^{\infty}_0 \mbb{E}_x[e^{-\lambda Z_{t}}1_{\{t<\zeta<\infty \}} ]dt\cr
\ar=\ar
\frac{1}{\mbb{P}_x(\zeta<\infty)}\int^{\infty}_0\mbb{E}_x[ e^{-\lambda Z_t}1_{\{t<\zeta\}}\mbb{E}_x [1_{\{\zeta<\infty\}}|\mcr{F}_t]] dt\cr
\ar=\ar
\frac{1}{\mbb{P}_x(\zeta<\infty)}\mbb{E}_x\Big[\int^{\zeta}_0 e^{-(\lambda+\Phi(0)) Z_{t}}dt\Big]\cr
\ar=\ar
\int^\infty_0 e^{-(\lambda+\Phi(0))y}[W(y)-e^{\Phi(0)x}W(y-x)]dy,
\eeqlb
where for the last equality we used (\ref{a00}) and (\ref{lap1}). Finally, thanks to (\ref{defscale})  we see that
(\ref{exp1}) holds.
\qed


\blemma\label{l2}
Let $f$ be a locally bounded function on $(0,\infty)$. Then for any $x > 0$  and $ \lambda> 0$, we have
$$
\mbb{E}_x\Big[\int^{\tau^-_y}_0 f(Z_{t}) dt\Big]<\infty\Leftrightarrow
\int^\infty f(z+y)[e^{-\Phi(0)(x-y)}W(z)-W(z-x+y)]dz<\infty.
$$
\elemma
\proof
 By the spatial homogeneousness of L\'{e}vy process and (\ref{lap1}) we have for $x>y>0$,
\beqlb\label{003}
\mbb{E}_x\Big[\int^{\tau^-_y}_0 f(Z_{t}) dt\Big]\ar=\ar
\mbb{E}_{x-y}\Big[\int^{\zeta}_0 f(Z_{t}+y) dt\Big]\cr
\ar=\ar
\int_0^\infty f(z+y)\mu_{x-y} (dz)\cr
\ar=\ar
\int_0^\infty f(z+y)[e^{-\Phi(0)(x-y)}W(z)-W(z-x+y)]dz.
\eeqlb
By the continuity of $W$, we have that function 
\[f(z+y)[e^{-\Phi(0)(x-y)}W(z)-W(z-x+y)]\]
 is bounded for $z$ near $0$. Then we only need to show the finiteness of the integral near $\infty$.
\qed

We are  now ready to prove the main results.

\noindent\emph{Proof of Theorem}~\ref{t1.1}.~
(i)$\Rightarrow$(iii)
For $x>0$ let $f_x(y):=f(y)1_{\{y\le x\}}$.  Given $x>0$, choose $\varepsilon>0$ small enough so that  $\varepsilon<x$.  From the non-negative jump property we see, for $t\leq \tau^-_\varepsilon$
$$f_x(Z_{t})\leq \sup_{y\geq \varepsilon}f(y)<\infty.$$ Therefore, by (\ref{exp1}) we have
\beqlb\label{infinity}
\mbb{E}_x\Big[\int^{\tau^-_\varepsilon}_0 f_x(Z_{t})e^{-\lambda Z_{t}}dt|\zeta<\infty\Big]
\ar\leq\ar
\sup_{y\geq \varepsilon}f(y)\mbb{E}_x\Big[\int^{\zeta}_0 e^{-\lambda Z_{t}}dt|\zeta<\infty\Big]\cr
\ar=\ar\sup_{y\geq \varepsilon}f(y)\frac{1-e^{-\lambda x}}{\psi(\lambda+\Phi(0))}
< \infty,
\eeqlb
If (i) holds, then there exists a $N>0$ such that
$$
\mbb{P}_x( A_{\zeta}(f_x)\le N|\zeta<\infty)>0.
$$
Define a stopping time
$$
T:=\inf\{t\geq 0: A_{t}(f_x)> N\}
$$
with the convention $\inf\emptyset=\infty.$
Then there exists a constant $\alpha\in (0,1]$ such that
\beqlb\label{0.4'}
\mbb{P}_x(T \ge \zeta|\zeta<\infty)=\alpha.
\eeqlb
Then we have for $0<\varepsilon<x$,
\beqlb\label{ineq1}\nonumber
\ar\ar \mbb{E}_x\Big[\int^{\tau^-_\varepsilon}_0 f_x(Z_{t})e^{-\lambda Z_{t}}dt| \zeta<\infty\Big]\\\nonumber
\ar\ar\qquad\qquad=
\mbb{E}_x\Big[\int^{\tau^-_\varepsilon}_0 f_x(Z_{t})e^{-\lambda Z_{t}}dt; T\geq \tau^-_\varepsilon| \zeta<\infty\Big]\\\nonumber
\ar\ar\qquad\qquad+
\mbb{E}_x\Big[\int^{T}_0 f_x(Z_{t})e^{-\lambda Z_{t}}dt;T< \tau^-_\varepsilon| \zeta<\infty\Big]\\\nonumber
\ar\ar\qquad\qquad+
\mbb{E}_x\Big[\int^{\tau^-_\varepsilon}_{T} f_x(Z_{t})e^{-\lambda Z_{t}}dt; T< \tau^-_\varepsilon| \zeta<\infty\Big]\\\nonumber
\ar\ar\qquad\qquad\leq
N\mbb{P}_x(T\geq \tau^-_\varepsilon| \zeta<\infty)\nonumber
+
N\mbb{P}_x(T< \tau^-_\varepsilon| \zeta<\infty)\\\nonumber
\ar\ar\qquad\qquad+
\mbb{E}_x\Big[\int^{\tau^-_\varepsilon}_{T} f_x(Z_{t})e^{-\lambda Z_{t}}dt;T< \tau^-_\varepsilon| \zeta<\infty\Big]\\
\ar\ar\qquad\qquad=
N+
\mbb{E}_x\Big[\int^{\tau^-_\varepsilon}_{T} f_x(Z_{t})e^{-\lambda Z_{t}}dt; T< \tau^-_\varepsilon| \zeta<\infty\Big].
\eeqlb
For each $\varepsilon<y\le x$,
\beqnn
\mbb{E}_x\Big[\int^{\tau^-_\varepsilon}_0 f_x(Z_{t})e^{-\lambda Z_{t}}dt| \zeta<\infty\Big]
\ar=\ar
 \mbb{E}_x\Big[\int^{\tau^-_y}_0 f_x(Z_t)e^{-\lambda Z_t}dt|\zeta<\infty\Big]\cr
 \ar\ar+
 \mbb{E}_x\Big[\int^{\tau^-_\varepsilon}_{\tau^-_y} f_x(Z_t)e^{-\lambda Z_{t}}dt| \zeta<\infty\Big].
\eeqnn
 By using Markov property, the second term of the above right hand side of the above equality can be written by
 \beqnn
  \mbb{E}_x\Big[\int^{\tau^-_\varepsilon}_{\tau^-_y} f_x(Z_t)e^{-\lambda Z_{t}}dt| \zeta<\infty\Big]
  \ar=\ar
   \frac{\mbb{E}_x\Big[\int^{\tau^-_\varepsilon}_{\tau^-_y} f_x(Z_t)e^{-\lambda Z_{t}}dt 1_{\{\zeta<\infty\}}\Big]}{\mbb{P}_y(\zeta<\infty)\mbb{P}_x(\tau_y^-< \infty)}\cr
   \ar=\ar
    \frac{\mbb{E}_x\left[\mbb{E}_x\Big[\int^{\tau^-_\varepsilon}_{\tau^-_y} f_x(Z_t)e^{-\lambda Z_{t}}dt 1_{\{\tau_y^-<\infty\}}1_{\{\zeta\circ\theta_{\tau^-_y}<\infty\}}|\mcr{F}_{\tau^-_y}\Big]\right]}{\mbb{P}_y(\zeta<\infty)\mbb{P}_x(\tau_y^-< \infty)}\cr
   \ar=\ar
    \mbb{E}_y\Big[\int^{\tau^-_\varepsilon}_{0} f_x(Z_{t})e^{-\lambda Z_{t}}dt| \zeta<\infty\Big],
 \eeqnn
 where $\theta$ is the shift operator. Then for each $\varepsilon<y<x$ we have
\beqlb\label{ineq3}
\mbb{E}_y\Big[\int^{\tau^-_\varepsilon}_0 f_x(Z_{t})e^{-\lambda Z_{t}}dt| \zeta<\infty\Big]\leq \mbb{E}_x\Big[\int^{\tau^-_\varepsilon}_0 f_x(Z_{t})e^{-\lambda Z_{t}}dt| \zeta<\infty\Big].
\eeqlb
Since $\varepsilon<Z_{T}<x$ on the event $\{T< \tau^-_\varepsilon\}$ under $\mbb{P}_x$, the strong Markov property and~(\ref{ineq3}) together yields
\beqlb\label{ineq2}
\nonumber\ar\ar \mbb{E}_x\Big[\int^{\tau^-_\varepsilon}_{T} f_x(Z_{t})e^{-\lambda Z_{t}}dt;  T< \tau^-_\varepsilon| \zeta<\infty\Big]\\
\ar\ar\qquad\qquad\nonumber
=\mbb{E}_x\Big[\int^{\tau^-_\varepsilon}_{T} f_x(Z_{t})e^{-\lambda Z_{t}}dt;  T< \tau^-_\varepsilon, \zeta<\infty\Big]/\mbb{P}_x(\zeta<\infty)\\
\ar\ar\qquad\qquad\nonumber
=\mbb{E}_x\Big\{1_{\{ T< \tau^-_\varepsilon\}}\mbb{E}_{Z_{T}}\Big[\int^{\tau^-_\varepsilon}_0 f_x(Z_{t})e^{-\lambda Z_{t}}dt; \zeta<\infty\Big]\Big\}/\mbb{P}_x(\zeta<\infty)\\
\ar\ar\qquad\qquad\nonumber
=\mbb{E}_x\Big\{1_{\{T< \tau^-_\varepsilon\}}\mbb{P}_{Z_T}(\zeta<\infty)\mbb{E}_{Z_T}\Big[\int^{\tau^-_\varepsilon}_0 f_x(Z_{t})e^{-\lambda Z_{t}}dt| \zeta<\infty\Big]\Big\}/\mbb{P}_x(\zeta<\infty)\\
\ar\ar\qquad\qquad
\leq
\mbb{P}_x(T< \tau^-_\varepsilon|\zeta<\infty)\mbb{E}_x\Big[\int^{\tau^-_\varepsilon}_0 f_x(Z_{t})e^{-\lambda Z_{t}}dt| \zeta<\infty\Big],
\eeqlb
where in the last inequality we used (\ref{ineq3}) and the fact that
\beqnn
\mbb{P}_x(T<\tau^-_\varepsilon, \zeta<\infty)
\ar=\ar
\mbb{E}_x\left[\mbb{E}_x[1_{\{T<\tau^-_\varepsilon, \zeta<\infty\}}|\mcr{F}_{T}]\right]\cr
\ar=\ar
\mbb{E}_x\left[1_{\{T<\tau^-_\varepsilon\}}\mbb{E}_x[1_{\{\zeta<\infty\}}|\mcr{F}_{T}]\right]\cr
\ar=\ar
\mbb{E}_x\left[1_{\{T<\tau^-_\varepsilon\}}\mbb{P}_{Z_T}(\zeta<\infty)\right].
\eeqnn
Combining ~(\ref{infinity}), (\ref{ineq1}) and~(\ref{ineq2})  we have
\beqnn
\ar\ar
\mbb{E}_x\Big[\int^{\tau^-_\varepsilon}_0 f_x(Z_{t})e^{-\lambda Z_{t}}dt|\zeta<\infty\Big]\\
\ar\ar\qquad\qquad\leq
N
+\mbb{P}_x(T< \tau^-_\varepsilon|\zeta<\infty)\mbb{E}_x\Big[\int^{\tau^-_\varepsilon}_0 f_x(Z_{t})e^{-\lambda Z_{t}}dt|\zeta<\infty\Big].
\eeqnn
It follows that
$$
\mbb{P}_x(T\geq \tau^-_\varepsilon|\zeta<\infty) \mbb{E}_x\Big[\int^{\tau^-_\varepsilon}_0 f_x(Z_{t})e^{-\lambda Z_{t}}dt\Big]
\leq
N.
$$
 Letting $\varepsilon\rightarrow0$ and using~(\ref{0.4'}) yields
$$
 \mbb{E}_x\Big[\int^{{\zeta}}_0 f_x(Z_{t})e^{-\lambda Z_{t}}dt|\zeta<\infty\Big]
\leq
N/\alpha<\infty.
$$
In addition,
\beqnn
\mbb{E}_x\Big[\int^{{\zeta}}_0 (f-f_x)(Z_{t})e^{-\lambda Z_{t}}dt|\zeta<\infty\Big]\leq \sup_{y\ge x}f(y)\mbb{E}_x\Big[\int^{\zeta}_0 e^{-\lambda Z_t}dt|\zeta<\infty\Big]<\infty,
\eeqnn
where the last inequality holds by Lemma \ref{l1'}.
Combing the above two inequalities gives (iii).

(iii)$\Rightarrow$(ii) From (iii) we have
$$\mbb{P}_x\Big[\int^{\zeta}_0 f(Z_{t})e^{-\lambda Z_{t}}dt<\infty|\zeta<\infty\Big]=1.$$
Since $\inf_{t\le \zeta}e^{-\lambda Z_{t}}>0 $ on event $\{\zeta<\infty\}$, we have a.s.
$$
\int^{\zeta}_0 f(Z_{t})dt<(\inf_{t\le \zeta}e^{-\lambda Z_{t}})^{-1}\int^{\zeta}_0 f(Z_{t})e^{-\lambda Z_{t}}<\infty
$$
on  event  $\{\zeta<\infty\}$. Then we obtain (ii).

(ii)$\Rightarrow$(i) Obvious.

(iii)$\Leftrightarrow$(iv) By (\ref{lap1}), we have
\beqlb\label{b1}
\ar\ar\mbb{E}_x\Big[\int^{\zeta}_0 f(Z_{t}) e^{-\lambda Z_{t}}dt|\zeta<\infty\Big]\cr
\ar=\ar
\mbb{E}_x\Big[\int^{\infty}_0 f(Z_{t}) e^{-\lambda Z_{t}}1_{\{t<\zeta<\infty\}}dt\Big]/\mbb{P}_x(\zeta<\infty)\cr
\ar=\ar
e^{\Phi(0)x}\int^{\infty}_0\mbb{E}_x\left[\mbb{E}_x[f(Z_{t}) e^{-\lambda Z_{t}}1_{\{t<\zeta<\infty\}}|\mcr{F}_t]\right]dt\cr
\ar=\ar
e^{\Phi(0)x}\int^{\infty}_0\mbb{E}_x\left[f(Z_{t}) e^{-\lambda Z_{t}}1_{\{t<\zeta\}}\mbb{E}_x[1_{\{\zeta<\infty\}}|\mcr{F}_t]\right]dt\cr
\ar=\ar
e^{\Phi(0)x}\int^{\infty}_0\mbb{E}_x[f(Z_{t}) e^{-(\lambda+\Phi(0)) Z_{t}}1_{\{t<\zeta\}}]dt\cr
\ar=\ar
e^{\Phi(0)x}\int^{\infty}_0f(y) e^{-(\lambda+\Phi(0))y}\mu_x(dy)\cr
\ar=\ar
\int^{\infty}_0f(y) e^{-(\lambda+\Phi(0))y}[W(y)-e^{\Phi(0)x}W(y-x)]dy.
\eeqlb
It follows from (\ref{defscale}) that
\beqnn
\int^{\infty}_0e^{-(\lambda+\Phi(0))y}[W(y)-e^{\Phi(0)x}W(y-x)]dy=\frac{1-e^{-\lambda x}}{\psi(\lambda+\Phi(0))}<\infty.
\eeqnn
Since for any $\varepsilon>0$, $f$ is bounded on $(\varepsilon,\infty)$ and $W(x)=0$ for $x<0$, the integral in (\ref{b1}) is finite if and only if
\beqnn
\int_{0+} f(y)W(y)dy<\infty.
\eeqnn
From (\ref{boundW}), the above inequality is equivalent to
$$\int_{0+} \frac{f(y)}{y\psi(1/y)}dy<\infty.$$
By a change of variable $y=1/\lambda$, we see the above inequality holds if and only if
\beqnn
\int^{\infty-}\frac{f(1/\lambda)}{\lambda\psi(\lambda)}d\lambda<\infty.
\eeqnn

\qed

\noindent\emph{Proof of Theorem}~\ref{t1.2}.~
(i)$\Rightarrow$(ii) We only need to prove that if $\mbb{P}_x(A_{\infty}(f)< \infty, \zeta= \infty)> 0$, then $\mbb{P}_x(A_{\infty}(f)= \infty, \zeta= \infty)= 0$.
If $\mbb{P}_x(A_{\infty}(f)< \infty, \zeta= \infty)> 0$, then there exists a $d> 0$ such that
$$\alpha:=\mbb{P}_x(A_{\infty}(f)\geq d)< 1, $$
where we make the convention that $f(x)=0$ for $x<0$. 
Let $\sigma_0:=0$ and for $n=0, 1, 2, \ldots$ define
$$\sigma_{n+1}:=\inf\left\{t\geq \sigma_n: A_t(f)> A_{\sigma_n}(f)+d, Z_{t\wedge\zeta}> x\right\}, \qquad n\geq 0.$$
 Since $Z_{\sigma_{k}} \geq x$ if $\sigma_k<\infty$, by the stochastic monotonicity of L\'{e}vy process and the monotonicity of $f$, we have for $\sigma_{k}<\infty$,
\begin{equation*}
\begin{split}
\mbb{P}_{Z_{\sigma_{k}}}(\sigma_{k+1} <\infty)
&\leq \mbb{P}_{Z_{\sigma_{k}}}(A_{\infty}(f)\ge d)\\
&\leq \mbb{P}_x (A_{\infty}(f)\ge d)=\alpha<1.
\end{split}
\end{equation*}
Notice  that $Z_t\to\infty$ if $\zeta=\infty$.  Using the strong Markov property and by induction, for any $n\ge 1$
we have
\beqnn
\mbb{P}_x(A_{\infty}(f)=\infty, \zeta= \infty)
 \ar\le\ar
\mbb{P}_x\bigg(\bigcap_{k=1}^n \big\{\sigma_k< \infty\big\}\bigg) \cr
 \ar=\ar
\mbb{E}_x\bigg[\prod_{k=1}^{n-1} 1_{\{\sigma_k< \infty\}} \mbb{P}_{Z_{\sigma_{n-1}}}\big(\sigma_n<\infty\big)\bigg] \cr
 \ar\le\ar
\alpha\mbb{E}_x\bigg[\prod_{k=1}^{n-1} 1_{\{\sigma_k< \infty\}}\bigg] \cr
 \ar= \ar
\alpha\mbb{P}_x\bigg(\bigcap_{k=1}^{n-1} \{\sigma_k< \infty\}\bigg)
 \le
\dots\le \alpha^n.
 \eeqnn
Letting $n\rightarrow\infty$, we have $\mbb{P}_x(A_{\infty}(f)=\infty, \zeta= \infty)=0.$

(ii)$\Rightarrow$ (iii)
 Fix a $x>0$ and $0<y<x$, by (ii) we have $\mbb{P}_y(A_{\infty}(f)<\infty|\zeta= \infty)=1$.
 Since $f$ is non-increasing, we can make the convention that $f(0)=\lim_{x\to0}f(x)>0$. Hence, if $\zeta<\infty$, then $A_{\infty}(f)=\infty$.
 Therefore,  there exists a constant $d>0$ and a constant $\alpha\in (0,1)$ such that
$$
\mbb{P}_y(A_{\infty}(f)<d)= \mbb{P}_y(A_{\infty}(f)<d, \zeta= \infty)>\alpha.
$$
Then by the stochastically monotonicity of L\'{e}vy process and the monotonicity of $f$ we have
$$
\mbb{P}_z(A_{\infty}(f)< d)\ge \mbb{P}_y(A_{\infty}(f)< d)>\alpha>0  \qquad\mbox{for each} ~z\ge y.
$$
It follows that
$$
\mbb{P}_z (A_{\infty}(f)\ge d)\leq 1-\alpha< 1   \qquad\mbox{for each} ~z\ge y.
$$
Then we have
$$
1>1-\alpha\ge \mbb{P}_z(A_{\infty}(f)\ge d)\ge \mbb{P}_z(A_{\tau^-_y}(f)\ge d)\qquad\mbox{for each} ~z\ge y.
$$
Put $\bar{\sigma}_0:=0$ and 
$$\bar{\sigma}_{n+1}:=\inf\left\{t> \bar{\sigma}_n: \int^{t}_{\bar{\sigma}_n}f(Z_{s})ds=d, Z_{s}> y~\mbox{for}~s\le t\right\}$$
for $n\geq 0$. 
Then we see that
$$
\{A_{\tau^-_y}(f)\ge nd\}=\{\bar{\sigma}_{n}<\infty\}.
$$
If for some $k\ge2$,  $\bar{\sigma}_{k-1}<\infty$, then we have $Z_{\bar{\sigma}_{k-1}}\ge y$ and
$$
\mbb{P}_{Z_{\bar{\sigma}_{k-1}}}(\bar{\sigma}_k< \infty)=\mbb{P}_{Z_{\bar{\sigma}_{k-1}}}(A_{\tau^-_y}(f)\ge d)\le 1-\alpha.
$$
Using the strong Markov property and combing the above formulas, by induction we have for $n\geq 1$,
\beqnn
\mbb{P}_x(A_{\tau^-_y}(f)\ge nd)
\ar=\ar
\mbb{P}_x\Big(\bigcap^n_{k=1} \{\bar{\sigma}_k< \infty\}\Big)\cr
\ar=\ar
\mbb{E}_x\Big[\mbb{P}_x\Big(\bigcap^n_{k=1} \{\bar{\sigma}_k< \infty\}\Big|\mcr{F}_{\bar{\sigma}_{n-1}}\Big)\Big]\cr
\ar=\ar
\mbb{E}_x\Big[\mbb{P}_{Z_{\bar{\sigma}_{n-1}}}(\bar{\sigma}_{n}< \infty)\prod^{n-1}_{k=1} 1_{\{\bar{\sigma}_k< \infty\}}\Big]\cr
\ar\le\ar
(1-\alpha)^n.
\eeqnn
It follows that $\mbb{E}_x [A_{\tau^-_y}(f)]<\infty$.

(iii) $\Rightarrow$ (i) Suppose that $\mbb{E}_x[A_{\tau^-_y}(f)]<\infty$ for all $x> y> 0$. Then $\mbb{P}_x(A_{\tau^-_y}(f)<\infty)=1.$  Since $\{\tau^-_y=\infty\}\subset \{\zeta=\infty\}$ and
\[\mbb{P}_x(\tau^-_y=\infty)=1-e^{-\Phi(0)(x-y)}>0\]
 for $x>y\ge 0$, we have
\[\mbb{P}_x(A_{\infty}(f)<\infty, \zeta=\infty)\ge \mbb{P}_x(A_{\infty}(f)<\infty, \tau^-_y=\infty)>0.\]
 It follows that $\mbb{P}_x(A_{{\infty}}(f)<\infty|\zeta= \infty)> 0.$

 Finally, under the condition $\psi'(0)<\infty$, we are going to show $(iii)\Leftrightarrow(iv)$.
 Thanks to Lemma~\ref{l2}, we only need to show that
 $$
\int^\infty f(z+y)[e^{-\Phi(0)(x-y)}W(z)-W(z-x+y)]dz<\infty\Leftrightarrow\int^\infty f(y)dy<\infty.
 $$
 By \cite[Lemma 2]{LiB19} we see that for any $x>y>0$,
 $$
 \lim_{z\to\infty} e^{-\Phi(0)(x-y)}W(z)-W(z-x+y)=\frac{1-e^{\Phi(0)(x-y)}}{\psi'(0)}.
 $$
Then we can immediately obtain the desired result.
\qed

\noindent\emph{Proof of Corollary}~\ref{cor2xx}.~
By (\ref{003}), we have
\beqnn
\ar\ar\int^{\infty}_0 f(z+y)[e^{-\Phi(0)(x-y)}W(z)-W(z-x+y)]dz\cr
\ar\ar\qquad=\int^{\infty}_0\int^\infty_0e^{-(z+y)\lambda}g(\lambda)[e^{-\Phi(0)(x-y)}W(z)-W(z-x+y)]d\lambda dz\cr
\ar\ar\qquad=\int_0^\infty e^{-\lambda y}g(\lambda)d\lambda\int^\infty_0 e^{-\lambda z} [e^{-\Phi(0)(x-y)}W(z)-W(z-x+y)]dz\\
\ar\ar\qquad=
\int_0^\infty e^{-\lambda y}\left(e^{-\Phi(0)(x-y)}-e^{-\lambda(x-y)}\right)\frac{g(\lambda)}{\psi(\lambda)}  d\lambda.
\eeqnn
Notice that $$\lim_{\lambda\to 0}e^{-\lambda y}\left(e^{-\Phi(0)(x-y)}-e^{-\lambda(x-y)}\right)= e^{-\Phi(0)(x-y)}-1$$ and $$e^{-\lambda y}\left(e^{-\Phi(0)(x-y)}-e^{-\lambda(x-y)}\right) \sim e^{-\lambda y}e^{-\Phi(0)(x-y)}\qquad \mbox{as} \qquad\lambda\to\infty.$$ Then we can finish the proof.\qed

\noindent\emph{Proof of Corollary}~\ref{Cor0}.~
(i) By letting $f(x)=x^{-\theta}$ in Theorem \ref{t1.1}, we immediately obtain the desired result.

(ii) For $f(x)=x^{-\theta}$, we have $f(x)=\frac{1}{\Gamma(\theta)}\int_0^\infty e^{-\lambda x} \lambda^{\theta-1}d\lambda.$ Then by letting $g(\lambda)=\frac{1}{\Gamma(\theta)}\lambda^{\theta-1}$ in Corollary \ref{cor2xx}, we  complete the proof.
 \qed

\end{document}